\documentclass[twocolumn,letterpaper,10pt]{article}
\usepackage{iasted}
\usepackage{times}
\usepackage[dvips]{graphicx}

\begin{document}
\parskip 0.05in
\parindent 0.2in 

\date{}

\title{Short-term equity dynamics and endogenous market fluctuations}

\author{
Ted V. Theodosopoulos \\
Department of Decision Sciences \\
and Department of Mathematics \\
Drexel University \\
Academic Building, Room 224 \\
Philadelphia, PA 19104, USA \\
email: theo@drexel.edu \\
\and
Muffasir H. Badshah \\
Department of Mathematics \\
Drexel University \\
email: mhb25@drexel.edu \\
}

\maketitle

\thispagestyle{empty}

\noindent
{\bf\normalsize ABSTRACT}\newline
{We present a model that investigates the spontaneous emergence of randomness in equity market microstructure. The phase space analysis of our model exposes an endogenous source of fluctuation in price and volume. We formulate a control problem for maximizing price regularity and stability while minimizing entanglement with the market, representing the NYSE specialists' affirmative obligation to maintain `fair and orderly markets'.} \vspace{2ex} 
   
\noindent
{\bf\normalsize KEY WORDS}\newline
{Endogenous uncertainty, equity market making}

\section{Overview}
Over the past ten years there has been an unprecedented increase in market-making activity, ranging from the globalized equity markets to progressively more customized structures that monetize novel risk dimensions.  Traditionally the domain of small, specialized, private partnerships, equity market-making in particular has attracted the attention of major financial intermediaries, and with it, the scrutiny of federal regulators and increasingly risk-conscious investors.

At the same time, a wealth of transaction-level data became widely accessible, enabling for the first time ever the pursuit of quantitative models for the short-term (high-frequency) dynamics of stock prices \cite{potters,bouchaud02,gourieroux,gu,hayashi}.  The interest in such models was further boosted by the apparent increase in the frequency and severity of well-publicized surprises that have rocked the global financial system recently, from the collapse of currencies to that of financial institutions and the bursting of equity and commodity bubbles.

Our interest is in the investigation of persistent dynamic effects that the rules of market making have on the prices of financial assets.  Also, we seek to construct a risk management framework applicable to equity market-making in general, and the monopolist NYSE specialists in particular.  Our contention is that the microstructure of the interactions between market-makers and investors represents an endogenous source of randomness \cite{shubik99} which percolates the financial intermediation network and gives rise to empirical mesoscopic phenomena that deviate from the predictions of the Efficient Markets Hypothesis (EMH).

The current paper represents a first step in the direction of modeling the short-term dynamics of equity market making.  Unlike earlier models that focus explicitly on the `microscopic' interactions between buyers and sellers \cite{bouchaud04}, our goal here is to attack the problem of the `ensemble' dynamics directly.  While we believe a complete theory will only emerge by exhibiting such ensemble dynamics as an appropriate `thermodynamic' (i.e. large scale) limit of some class of microscopic interaction models \cite{day91}, this goal appears still out of reach.  We see our work presented here as an attempt to `bridge' the gap between the burgeoning literature on constructive microscopic models and the increasingly refined empirical studies of `stylized facts' about the equity markets' deviations from the EMH \cite{gabaix}.

We begin by presenting our model as a coupled system of nonlinear ordinary differential equations (ODEs).  The state of our system encodes the intra-day price variation of an equity, the innovations in its instantaneous trading volume, and two slope parameters that couple dynamically the fluctuations of price and volume.  We proceed to analyze this model using techniques from dynamical systems theory.  Our main result is that, in general, bifurcations in phase space lead to endogenous uncertainty in the simultaneous determination of price and volume.  We proceed to formulate a decision problem faced by the specialist as they fulfill their affirmative obligation to the NYSE to maintain price regularity while minimizing their interference with the market.  We end by presenting our next goals in this research program.

\section{Model}
The model we present here was motivated by a series of papers by Day and collaborators that investigate plausible mechanisms that lead to endogenous fluctuations in macroeconomics \cite{day85,day87,saari} and financial transaction data \cite{day93,gu,shubik97}.  We proceed by constructing the following set of coupled differential equations for the tick-by-tick price process $x_1$ and the daily trading volume $x_2$:
\begin{eqnarray}
\dot{x}_1 & = & x_4^{-1} x_1 + z \label{eq:xdot} \\
\dot{x}_2 & = & x_3 x_4^{-1} x_1 \label{eq:ydot} \\
\dot{x}_3 & = & \beta_1 x_2 \label{eq:alphadot} \\
\dot{x}_4 & = & \beta_2 x_2 \label{eq:gammadot}
\end{eqnarray}
where $x_3$ and $x_4$ are unobservable components of the state, and $z$ denotes an exogenous input, which can be used to model the specialist's control the price setting process.  

Equation (\ref{eq:xdot}) provides a first-order response dynamics for price momentum, which is justified for short enough time windows as the result of a Taylor expansion \cite{bouchaud04}.  Specifically, one can think of (\ref{eq:xdot}) as capturing a classical `momentum market' when $x_4 >0$ and a `contrarian market' when $x_4 <0$.  Equation (\ref{eq:ydot}) captures the first order condition around the local supply/demand equilibrium\footnote{Observe that, together with (\ref{eq:xdot}), this implies that $\dot{x}_2 = x_3 \dot{x}_1$.}.  In this context, $x_3<0$ generally signifies a supply shock, while $x_3>0$ signifies a demand shock.  Equation (\ref{eq:alphadot}) captures the second order condition, with 
$${\frac {d^2 x_1}{d x_2^2}} = -{\frac {\beta_1 x_2}{x_3^2 \dot{x}_2}}.$$
This implies for instance that $\beta_1>0$ means a supply shock has larger price effects but smaller volume effects than a demand shock.  On the other hand equation (\ref{eq:gammadot}) specifies the convexity of the phase picture for the price process $x_1$.

\section{Master Equation}

Dividing (\ref{eq:alphadot}) and (\ref{eq:gammadot}) and integrating we obtain
\begin{equation}
{\frac {x_3}{x_4}} = {\frac {\beta_1}{\beta_2}} + {\frac {x_3(0) - {\frac {\beta_1}{\beta_2}} x_4(0)}{x_4}}.  \label{eq:x3overx4}
\end{equation}
From now on we will assume that the initial conditions for (\ref{eq:xdot}-\ref{eq:gammadot}) are such that $\beta_3 x_3(0) = \beta_1 x_4(0)$.  Then (\ref{eq:x3overx4}) simplifies to $x_3/x_4 = \beta_1/\beta_2$, transforming (\ref{eq:ydot}) to $\ddot{x}_2 = \dot{x}_1 \beta_1 \beta_2^{-1}$.
On the other hand, we see that 
\begin{eqnarray}
\ddot{x}_4 & = & \beta_2 \dot{x}_2  \nonumber \\
& = & \beta_1 \dot{x}_1 x_4 - \beta_1 x_4 z  \nonumber \\
& = & \beta_2 \ddot{x}_2 x_4 - \beta_1 x_4 z  \nonumber \\
& = & x_4 \left( \dot{\ddot{x}}_4 - \beta_1 z \right)  \label{eq:master1}
\end{eqnarray}
where the first equality follows by differentiating (\ref{eq:gammadot}) and the final one by differentiating (\ref{eq:gammadot}) once more.  Focusing on the homogeneous equation ($z \equiv 0$), integration by parts of (\ref{eq:master1}) leads to an integral of motion:
\begin{equation}
2 x_4 \ddot{x}_4 - \dot{x}_4 \left( \dot{x}_4 + 2\right) = C_1  \label{eq:master2}
\end{equation}
where $C_1$ is the initial value of the left-hand-side.  Thus, the operator ${\cal L}: f \mapsto 2f \ddot{f} - \dot{f} \left( \dot{f} +2 \right)$ represents a (nonlinear) symmetry of the short-term equity dynamics desribed by (\ref{eq:xdot}-\ref{eq:gammadot}).  Observe that this one-dimensional second-order nonlinear ODE applies to any stock, irrespective of its type.  The parameters $\beta_1$ and $\beta_2$ enter in the computation of $x_1$ and $x_2$ from $x_4$ and its derivatives:
\begin{eqnarray}
x_1 & = & \beta_1^{-1} \ddot{x}_4 \\  \label{eq:xfromgammaddot}
x_2 & = & \beta_2^{-1} \dot{x}_4  \label{eq:yfromgammadot}
\end{eqnarray}
Letting $u=x_4$ and $v = \dot{u} = \dot{x}_4$, we can transform (\ref{eq:master2}) into a separable ODE by observing that 
$${\frac {dv}{du}} = {\frac {\ddot{x}_4}{\dot{x}_4}} \Rightarrow \ddot{x}_4 = v {\frac {dv}{du}},$$
which, when substituted into (\ref{eq:master2}) yields
\begin{equation}
{\frac {dv}{du}} = {\frac {C_1+ v(v+2)}{2uv}}  \label{eq:master3}
\end{equation}

\section{Phase Picture}

Integrating (\ref{eq:master3}) with respect to $u$ and assuming $v \left(u=1 \right) = v_0$ we obtain
$$\int_{v_0}^v {\frac {2t dt}{C_1 + t(t+2)}} = \int_1^u {\frac {du}{u}} \Longrightarrow$$
$$
\Longrightarrow \log \left({\frac {v^2 +2v +C_1}{v_0^2 +2v_0 +C_1}} \right) - 2\int_{v_0}^v {\frac {dt}{t^2 + 2t +C_1}} = \log u \Longrightarrow$$
$$\Longrightarrow |u| = \left\{
\begin{array}{lll}
C_2 \left( v^2 +2v +C_1 \right) \exp \left\{ -h(x)\right\} & \mbox{if $C_1>1$} \\
C_4 (v + 1)^2 \exp \left\{ - {\frac {2v}{v +1}} \right\} & \mbox{if $C_1=1$} \\
g(v) & \mbox{otherwise}
\end{array}
\right.,$$
where
$$h(x)=   {\frac {2}{\sqrt{C_1 -1}}} \arctan \left( {\frac {v+1}{\sqrt{ C_1 -1}}} \right),$$
$$C_5=\left(1-C_1 \right)^{-1/2} \mbox{   when $C_1 \leq 1$},$$
$$\ell_+ (v) = \left| v+1+C_5^{-1} \right|,$$
$$\ell_- (v) = \left| v+1-C_5^{-1} \right| \mbox{    and}$$
$$g(v) = \left\{
\begin{array}{ll}
C_3 \ell_{\scriptscriptstyle +}^{\scriptscriptstyle 1+ C_5} \ell_{\scriptscriptstyle -}^{\scriptscriptstyle -1+ C_5} & \mbox{if $C_1<1$ and $v \not\in A$} \\
-C_3 \ell_{\scriptscriptstyle +}^{\scriptscriptstyle 1+ C_5} \ell_{\scriptscriptstyle -}^{\scriptscriptstyle -1+ C_5} & \mbox{if $C_1<1$ and $v \in A$}
\end{array}
\right.,$$
$A = \left. \left[-1 - \sqrt{1- C_1}, -1 + \sqrt{1- C_1} \right. \right)$, and $C_2$, $C_3$ and $C_4$ are the natural normalization constants.

\begin{figure}
\centering
\includegraphics[width=2.5in]{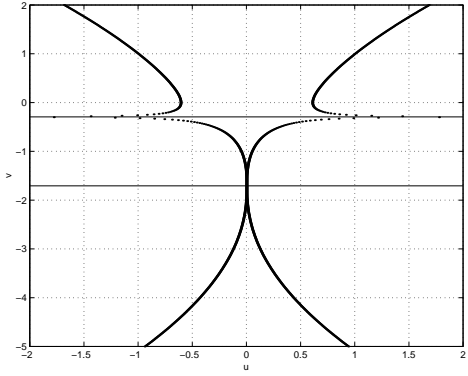}
\caption{Phase diagram for $C_1=0.5$ and $v_0=1$.  The two horizontal lines correspond to $-1 \pm \sqrt{2}/2$.}
\label{fig:plushalf}
\end{figure}

\begin{figure}
\centering
\includegraphics[width=2.5in]{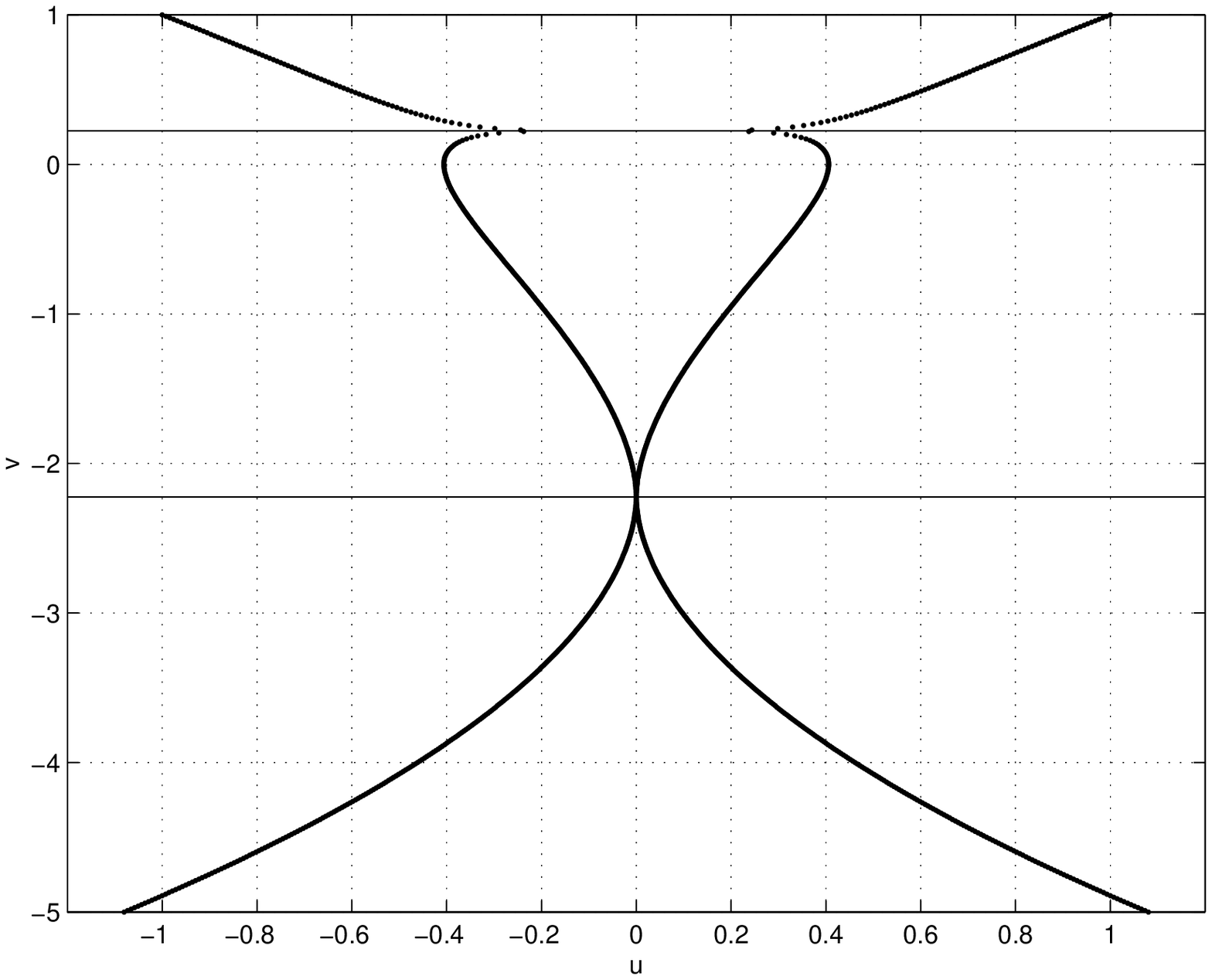}
\caption{Phase diagram for $C_1=-0.5$ and $v_0=1$.  The two horizontal lines correspond to $-1 \pm \sqrt{6}/2$.}
\label{fig:minushalf}
\end{figure}

\begin{figure}
\centering
\includegraphics[width=2.5in]{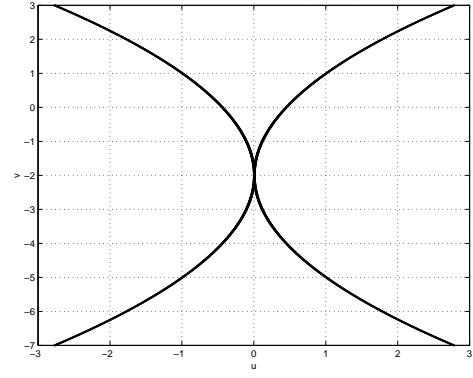}
\caption{Phase diagram for $C_1=0$ and $v_0=1$.}
\label{fig:zero}
\end{figure}

\begin{figure}
\centering
\includegraphics[width=2.5in]{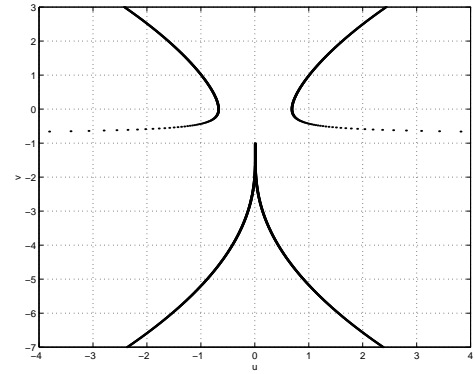}
\caption{Phase diagram for $C_1=1$ and $v_0=1$.}
\label{fig:one}
\end{figure}

\begin{figure}
\centering
\includegraphics[width=2.5in]{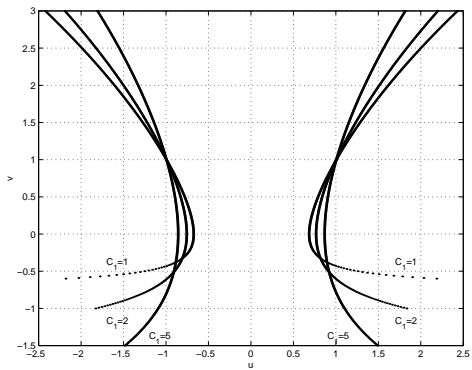}
\caption{Phase diagrams for $C_1=1,2,5$ and $v_0=1$}
\label{fig:morethanone}
\end{figure}

From the above analysis we can infer the following properties of the resulting phase picture:

\begin{enumerate}
\item $\lim_{|u| \rightarrow \infty} |v| = \infty$.
\item When $C_1 > 1$, the phase diagram consists of two disconnected symmetric components on either side of the $u=0$ axis (see Figure \ref{fig:morethanone}).
\item When $C_1 \neq 0$, the phase diagram crosses the $v=0$ axis perpendicularly, i.e. $\left| {\frac {dv}{du}}\right|_{v=0} = \infty$ (see Figures \ref{fig:plushalf}, \ref{fig:minushalf} and \ref{fig:morethanone}).  On the other hand when $C_1=0$, $\left| {\frac {dv}{du}}\right|_{v=0} = \pm u$ (see Figure \ref{fig:zero}).
\item When $C_1 \in (-\infty, 1) \setminus \{0\}$, the phase diagram has a unique pole at $-1+\sqrt{1-C_1}$ (see Figures \ref{fig:minushalf}, \ref{fig:plushalf} and \ref{fig:one}).
\item When $C_1 \leq 1$, the phase diagram is tangent to the $u=0$ axis at $-1-\sqrt{1-C_1}$ (see Figures \ref{fig:plushalf}, \ref{fig:minushalf}, \ref{fig:zero} and \ref{fig:one}).
\end{enumerate}

\vspace{-0.1in}

We will use Figure \ref{fig:minushalf} to exhibit the possible solution paths.  Specifically, we interpret $v$ as the time derivative of $u$ and follow the dynamics:

\begin{enumerate}
\item When $u>0$ and $v>-1+\sqrt{6}/2$, $u$ explodes to $+\infty$.  Similarly, when $u<0$ and $v<-1-\sqrt{6}/2$, $u$ explodes to $-\infty$.
\vspace{-0.1in}
\item Let $B= \left( 0,-1-\sqrt{6}/2 \right)$, $C=(-0.4,0)$, $D=(0.4,0)$ and $E= \left(0, -1+\sqrt{6}/2 \right)$.  These are bifurcation points, in that they entail choices for the path.  Specifically, points $B$ and $E$ are sensitive to small perturbations along the $v$ axis.  Using (\ref{eq:xfromgammaddot}) we see that these perturbations measure the value of $x_1$.  On the other hand, points $C$ and $D$ are sensitive to small perturbations along the $u$ axis.  Using (\ref{eq:yfromgammadot}) we see that these perturbations measure the value of $x_2$.
\vspace{-0.1in}
\item Let $F$ denote the point on the phase diagram in the fourth quadrant with $u=0.4$.  Also, let $G$ denote the point on the phase diagram in the second quadrant with $u=-0.4$.  Four cycles are possible in this phase diagram, depending on choice made at the $x_2$-bifurcation points: $BCEDB$, $BCGEDB$, $BCEDFB$ and $BCGEDFB$.  On the other hand the $x_1$-bifurcation points can cause explosions, depending on the signs of $\beta_1$ and $\beta_2$.
\end{enumerate}
\vspace{-0.1in}

\begin{figure}
\centering
\includegraphics[width=2.5in]{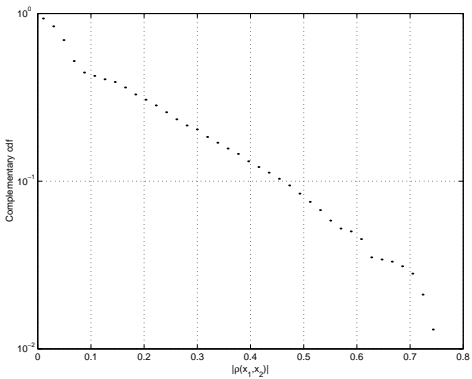}
\caption{Histogram of absolute value of the correlation between $x_1$ and $x_2$ from a Monte Carlo simulation of the dynamics shown in Figure \ref{fig:minushalf} using 1,000 randomly selected initial conditions.}
\label{fig:correl}
\end{figure}

All the versions of the phase diagram (i.e. for any value of $C_1$) exhibit similar bifurcation properties.  The resulting sensitivity to small perturbations characterizes the system dynamics.  Furthermore, the resulting path dependence on the value of $x_1$ and $x_2$ renders it impossible to associate a unique value to one if we know the value of the other.  The residual randomness that is accumulated through the time evolution of this system is what we referred to earlier as an uncertainty principle.  Figure \ref{fig:correl} shows that the absolute value of the correlation $\rho \left(x_1, x_2 \right)$ for $C_1=-0.5$ and randomly chosen starting configurations decreases exponentially at first, and stays strictly bounded away from 1 (in fact, lower than $0.8$).  The correlation distributions for different values of $C_1$ are qualitatively similar.

\section{Volume and Price Distributions}

\begin{figure}
\centering
\includegraphics[width=2.5in]{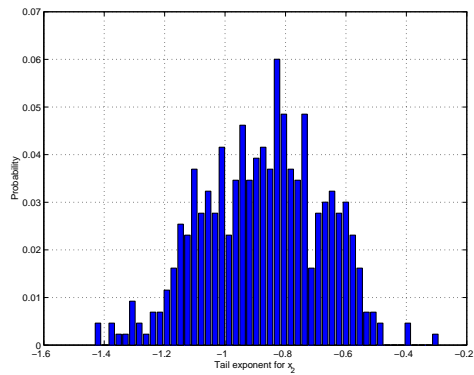}
\caption{Histogram of $\lambda_2$ from a Monte Carlo simulation of the dynamics shown in Figure \ref{fig:minushalf} using 1,000 randomly selected initial conditions.}
\label{fig:volumehist}
\end{figure}

\begin{figure}
\centering
\includegraphics[width=2.5in]{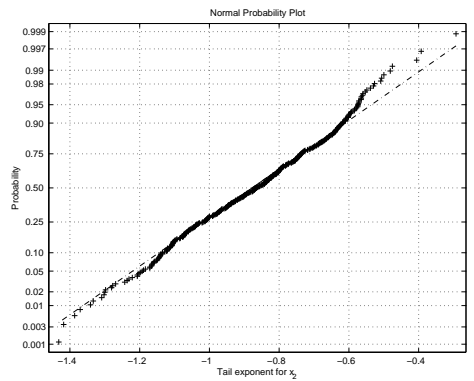}
\caption{Normal probability plot of $\lambda_2$ from a Monte Carlo simulation of the dynamics shown in Figure \ref{fig:minushalf} using 1,000 randomly selected initial conditions.}
\label{fig:volumenorm}
\end{figure}

We have used the Monte Carlo simulation for the dynamics in the case $C_1=-0.5$, $\beta_1= \beta_2 =1$ to compute the distribution of volume and price.  The tails of each of the resulting distributions was fit to a power law and the corresponding exponents were recorded.  Specifically, for each Monte Carlo run $i$, we computed the empirical measure of $x_2(t)$ and $x_1(t)$ and computed $\lambda_1(i)$ and $\lambda_2(i)$ so that $\Pr_i \left(x_j < a \right) \sim a^{\lambda_j(i)}$ for $j=1,2$.  Figure \ref{fig:volumehist} shows the distribution of the resulting exponents $\lambda_2 (i)$ for the left tail of the volume distribution.  The exponents for the right tail of the volume distribution exhibit similar qualitative characteristics (i.e. normality, but with different mean and variance).
We see that the distribution of $\lambda_2$ is approximately normal.  Indeed, Figure \ref{fig:volumenorm} shows the corresponding normal probability plot, which confirms our suspicion that the distribution of volume exponents is approximately normal, with mean $\mu_2 \cong -0.88$ and standard deviation $\sigma_2 \cong 0.12$.  
On the other hand the distribution of the tail exponents for prices has a pronounced negative skew, with a mean $\mu_1 \cong -0.56$ and standard deviation $\sigma_1 \cong 0.33$.
\begin{figure}
\centering
\includegraphics[width=2.5in]{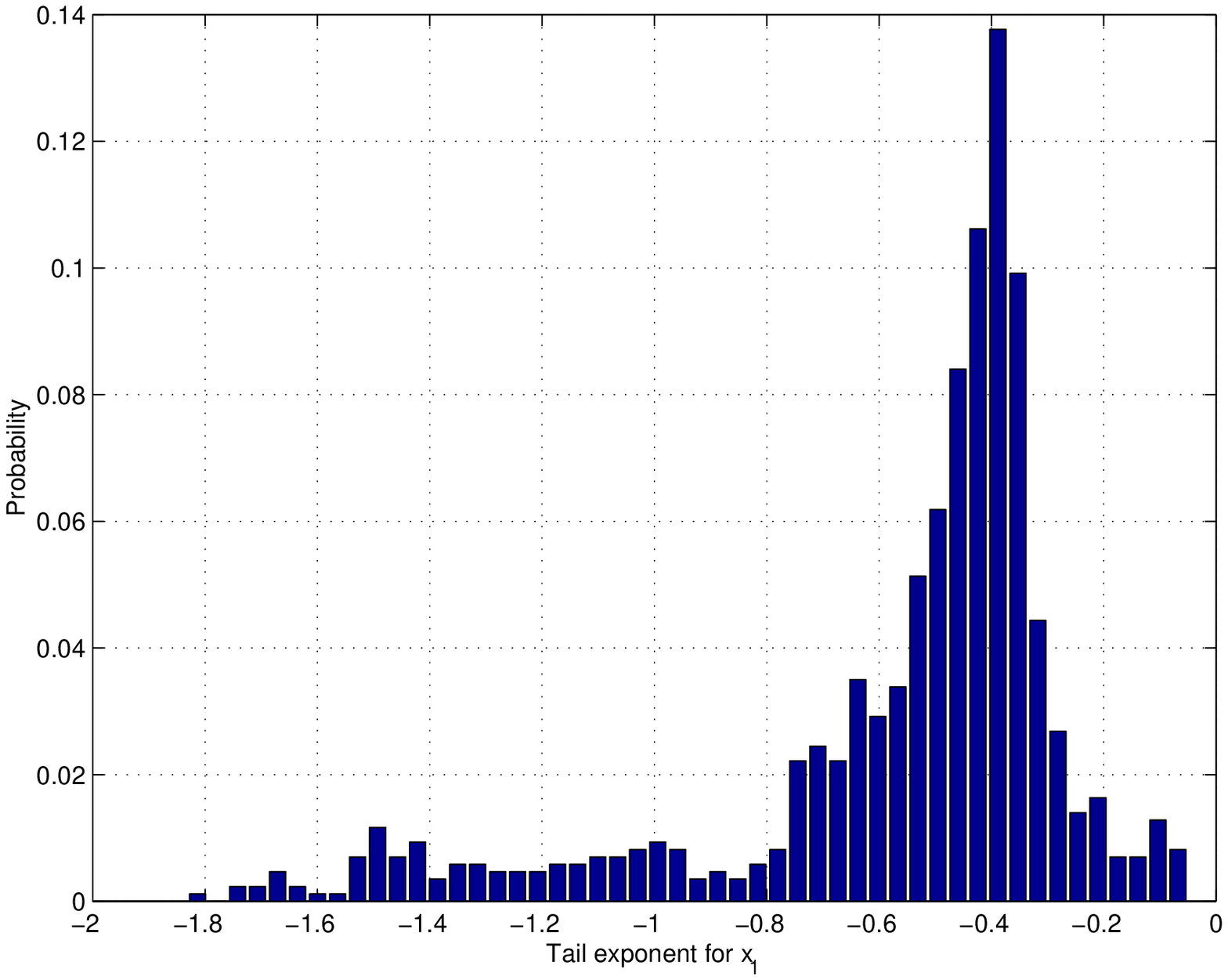}
\caption{Histogram of $\lambda_1$ from a Monte Carlo simulation of the dynamics shown in Figure \ref{fig:minushalf} using 1,000 randomly selected initial conditions.}
\label{fig:pricehist}
\end{figure}
In both cases, the absolute values of the exponents are substantially less than 2, implying light-tailed distributions, both for the volume and the price returns.  However, these distributions are not directly comparable to empirically determined distributions.  It should be remembered that the randomness we observe here arises entirely endogenously, through the dynamic bifurcations in phase space, rather than through any exogenous inputs.  Figures \ref{fig:volumehist} and \ref{fig:pricehist} are shown as examples.  Different values of $C_1$, $\beta_1$ and $\beta_2$ lead to different power distributions, with similar qualitative characteristics.  
\vspace{-0.1in}
\section{Exogenous Input}

Observe that for any exogenous input $z$, $\int u(t) z(t) dt = \int u v^{-1} \tilde{z} (u) du$ using the fact that $du=v dt$, where $u(t) = u \Rightarrow \tilde{z} (u) = z(t)$.  If we now return to (\ref{eq:master2}) and add the effect of the exogenous input $z$ from (\ref{eq:master1}) we obtain in lieu of (\ref{eq:master3}) the following equation:
\begin{equation}
{\frac {dv}{du}} = {\frac {C_1+ v(v+2) + 2 \beta_1 \int u v^{-1} \tilde{z} du}{2uv}}  \label{eq:master4}
\end{equation}
Differentiating with respect to $u$ we obtain after some algebra:
\begin{equation}
uv^2 {\frac {d^2 v}{du^2}} + uv \left( {\frac {dv}{du}} \right)^2 – v {\frac {dv}{du}} -\beta_1 u \tilde{z} = 0  \label{eq:master5}
\end{equation}
Equation (\ref{eq:master5}) provides us with a tool to choose a control input $z$ as a function of $u$ in order to attain a desirable $(u,v)$ profile.  Of course, if we choose $\tilde{z} (u) \equiv 0$, after some algebra we recover (\ref{eq:master1}).  On the other hand, if we choose for instance $\tilde{z} (u) = ku^{-1}$, then (\ref{eq:master5}) simplifies to the following system of two coupled first order ODEs:
\begin{eqnarray*}
& & {\frac {dw}{du}} = {\frac {vw(1 - uw) + k\beta_1}{u v^2}} \\
& & {\frac {dv}{du}} = w,
\end{eqnarray*}
which can be integrated numerically.

We are finally in a position to formulate the control problem \cite{runggaldier} faced by a monopolist market maker at the NYSE.  In particular, the specialists operate under an `affirmative' obligation by the exchange to maintain fair and orderly markets.  Generally, this requirement is interpreted as an obligation to complete any one-sided market and minimizing price variation while avoiding as much as possible interference with the autonomous evolution of the market.

The model we developed in this paper allows us to formalize this decision problem.  The price regularity constraint can be conceptualized as an upper bound $U$ to the integrated square increments of $x_1$, $\int \left(\dot{x}_1 \right)^2 dt$.  This constraint is often quantified by the exchange as the maximum allowable bid-offer spread, which can vary across stocks and from time to time, depending on the observed trading patterns.

The specialist's objective as a market maker under the NYSE affirmative obligation is the minimization of the `entanglement' with the market, which can be measured as the absolute value of the correlation between the specialist's actions $z$ and the resulting price series $x_1$ (see also the discussion in \cite{day93,gu}).  On the other hand the specialist's profit, $\Pi$, is generally a multiple of the trading volume, $x_2$, so long as the price does not exhibit a pronounced trend.  A trending market is detrimental to the specialist's profit because it often necessitates an inordinate number of one-sided trades against the market.  Thus, we can model the rate of change of $\Pi$ as
\begin{eqnarray}
\dot{\Pi} & = & c_1 U \dot{x}_2 - c_2 \dot{x}_1 \nonumber \\
& = & c_1 U \beta_1 \beta_2^{-1} x_1 - c_2 \dot{x}_1,	\label{eq:profit}
\end{eqnarray}
for some positive constants $c_1, c_2>0$.
In summary, we propose the following control problem for the monopolist market maker:
\begin{eqnarray}
& \min_{\displaystyle z(\cdot)} & \rho(z,x_1) \nonumber \\
& \mbox{subject to} & \nonumber \\
& & uv^2 {\frac {d^2 v}{du^2}} + uv \left( {\frac {dv}{du}} \right)^2 - v{\frac {dv}{du}} -\beta_1 u \tilde{z} = 0 \nonumber \\
& & x_1 = \beta_1^{-1} {\frac {dv}{du}} \nonumber \\
& & \int \left(\dot{x}_1 \right)^2 dt \leq U \nonumber \\
& & c_2 \beta_2 \dot{x}_1 - c_1 \beta_1 U x_1 + \beta_2 L \leq 0,  \label{eq:control}
\end{eqnarray}
where $L$ is a lower bound to $\dot{\Pi}$ from (\ref{eq:profit}) that the specialist is willing to accept (loss rates higher than $|L|$ are unacceptable and trigger restructurings that are outside the scope of the current model).
\vspace{-0.1in}
\section{Conclusions and Next Steps}

We presented a dynamical model (\ref{eq:xdot}-\ref{eq:gammadot}) for the short-term behavior of stocks.  The formal analysis of the model led to the determination of an integral of motion, $C_1$, which determines the qualitative behavior of the dynamics, as shown in Figures \ref{fig:plushalf}, \ref{fig:minushalf}, \ref{fig:zero}, \ref{fig:one} and \ref{fig:morethanone}.

Investigations of the homogeneous phase picture obtained in (\ref{eq:master3}) helped us identify bifurcation points, which introduce an endogenous source of uncertainty in the high-frequency dynamics of price and volume \cite{day93,shubik99}.  Using this framework, we formulated a control problem for the monopolist market maker, that attempts to balance the desire for price regularity and profit making with the obligation to avoid interference with the market (\ref{eq:control}).

The solution of the specialist's control problem is a natural next step of the analysis presented here.  This problem can be seen from the perspective of the market maker, as presented above, in which case we are seeking an optimal `market making strategy', $z$.  Alternatively, this problem can be seen from the perspective of a risk manager responsible for the capital adequacy of the market maker.  In this case it is the maximum allowable loss rate $L$ that is of primary interest.  Yet a different perspective on this problem is that of the regulators who design rules for the exchange.  In this view, it is the price regularity measure $U$ that is central, as well as the minimum achievable market interference (optimal value of the objective function).

In a different direction, it is desirable to investigate how the endogenous randomness introduced by the nonlinear dynamics (\ref{eq:master2}) interacts with the exogenous input $z$, and how it scales with time.  To begin with we may consider the dynamics presented here as a nonlinear filter, taking as an input random paths $z$ and producing filtered versions.  The study of this filtering paradigm entails the investigation of (\ref{eq:master5}).  Further, even when $z$ is not random but the result of an optimization as shown in (\ref{eq:control}), it will still interact with the randomness produced by the homogeneous dynamics, because $z$ effectively continuously readjusts $C_1$ as we can see by comparing (\ref{eq:master4}) and (\ref{eq:master3}).  Finally one would wish to recover in some appropriate `long time' limit the tail exponents for the volume and returns distributions that are empirically observed \cite{gabaix,hayashi,bouchaud04}.

\end{document}